\documentclass[11pt, reqno]{amsart}

\usepackage{amssymb}
\usepackage[T1]{fontenc}
\usepackage{pdfpages}

\evensidemargin\oddsidemargin

\begin{document}

\baselineskip=18pt
\setcounter{page}{1}
    
\newtheorem{Conj}{Conjecture\!\!}
\newtheorem{TheoA}{Theorem A\!\!}
\newtheorem{TheoB}{Theorem B\!\!}
\newtheorem{TheoC}{Theorem C\!\!}
\newtheorem{TheoD}{Theorem D\!\!}
\newtheorem{Lemm}{Lemma}
\newtheorem{Rem}{Remark}
\newtheorem{Coro}{Corollary}
\newtheorem{Propo}{Proposition}

\renewcommand{\theConj}{}
\renewcommand{\theTheoA}{}
\renewcommand{\theTheoB}{}
\renewcommand{\theTheoC}{}
\renewcommand{\theTheoD}{}

\def\a{\alpha}
\def\b{\beta}
\def\B{{\bf B}} 
\def\CC{{\mathbb{C}}} 
\def\cN{{\mathcal{N}}} 
\def\cH{{\mathcal{H}}} 
\def\cI{{\mathcal{I}}} 
\def\cS{{\mathcal{S}}}
\def\UU{{\mathcal{U}}}
\def\ca{c_{\a}}
\def\ka{\kappa_{\a}}
\def\coa{c_{\a, 0}}
\def\cua{c_{\a, u}}
\def\cL{{\mathcal{L}}} 
\def\cW{{\mathcal{W}}} 
\def\Ea{E_\a}
\def\Eab{E_{\a,\b}}
\def\eps{{\varepsilon}} 
\def\esp{{\mathbb{E}}} 
\def\Ga{{\Gamma}} 
\def\G{{\bf \Gamma}} 
\def\K{{\bf K}}
\def\HH{{\bf H}}
\def\ii{{\rm i}}
\def\e{{\rm e}}
\def\L{{\bf L}}
\def\lbd{\lambda}
\def\lacc{\left\{}
\def\lcr{\left[}
\def\lpa{\left(}
\def\lva{\left|}
\def\M{{\bf M}}
\def\Ma{\M_\a}
\def\Mab{\M_{\a,\b}}
\def\Mabe{\M_{\a,\b, \varepsilon}}
\def\NN{{\mathbb{N}}} 
\def\pb{{\mathbb{P}}} 
\def\tpab{\hat{\varphi}_{a,b}} 
\def\tpa{\hat{\psi}_{\a}}
\def\tppa{\tilde{\psi}_{\a}} 
\def\tva{\hat{\varphi}_{\a}} 
\def\rl{{\mathbb{R}}}
\def\racc{\right\}}
\def\rpa{\right)}
\def\rcr{\right]}
\def\rva{\right|}
\def\prost{{\succ_{\! st}}}
\def\W{{\bf W}}
\def\X{{\bf X}}
\def\Y{{\bf Y}}
\def\Z{{\bf Z}}
\def\Za{\Z_{\a}}
\def\Xab{\X_{\a,\b}}
\def\Yab{\Y_{\a,\b}}
\def\Yabe{\Y_{\a,\b, \varepsilon}}
\def\XX{{\mathcal X}}
\def\Y{{\bf Y}}
\def\U{{\bf U}}
\def\V{{\bf V}_\a}
\def\Un{{\bf 1}}
\def\ZZ{{\mathbb{Z}}}
\def\A{{\bf A}}
\def\AA{{\mathcal A}}
\def\hAA{{\hat \AA}}
\def\hL{{\hat L}}
\def\hT{{\hat T}}

\def\claw{\stackrel{d}{\longrightarrow}}
\def\elaw{\stackrel{d}{=}}
\def\qed{\hfill$\square$}

\newcommand*\pFqskip{8mu}
\catcode`,\active
\newcommand*\pFq{\begingroup
        \catcode`\,\active
        \def ,{\mskip\pFqskip\relax}%
        \dopFq
}
\catcode`\,12
\def\dopFq#1#2#3#4#5{%
        {}_{#1}F_{#2}\biggl[\genfrac..{0pt}{}{#3}{#4};#5\biggr]%
        \endgroup
}

\title{On the log-concavity of the Wright function}

\author[Rui. A. C. Ferreira]{Rui A. C. Ferreira}

\address{Grupo F\'isica-Matem\'atica, Faculdade de Ci\^encias, Universidade de Lisboa,
Av. Prof. Gama Pinto, 2, 1649-003, Lisboa. {\em Email}: {\tt raferreira@fc.ul.pt}}

\author[Thomas Simon]{Thomas Simon}

\address{Laboratoire Paul Painlev\'e, Universit\'e de Lille, 59000 Lille. {\em Email}: {\tt thomas.simon@univ-lille.fr}}



\begin{abstract} 
We investigate the log-concavity on the half-line of the Wright function $\phi(-\a,\b,-x),$ in the probabilistic setting $\a\in (0,1)$ and $\b \ge 0.$ Applications are given to the construction of generalized entropies associated to the corresponding Mittag-Leffler function. A natural conjecture for the equivalence between the log-concavity of the Wright function and the existence of such generalized entropies is formulated. The problem is solved for $\b\ge\a$ and in the classical case $\b = 1-\a$ of the Mittag-Leffler distribution, which exhibits a certain critical parameter $\a_*= 0.771667...$ defined implicitly on the Gamma function and characterizing the log-concavity. We also prove that the probabilistic Wright functions are always unimodal, and that they are multiplicatively strongly unimodal if and only if $\b\ge\a$ or $\a\le 1/2$ and $\b = 0.$ 
\end{abstract}

\maketitle

\section{Introduction}

The Wright function is defined as the entire function  
$$\phi (\rho, \beta, z) \; =\; \sum_{n\ge 0} \frac{z^n}{n!\, \Gamma (\beta + \rho n)}, \qquad  \beta, z\in \CC,\; \rho > - 1.$$
In the case $\rho = 0,$ it is an exponential function and in the other cases it is an entire function of order $1/(1+\rho).$ In the early literature, the Wright function is sometimes called the generalized Bessel function because of the following representation in the case $\rho = 1$:
$$J_\nu(z) \; =\; (z/2)^\nu \phi (1, \nu + 1, - z^2/4).$$ The Wright function admits an integral representation given by
\begin{equation}
\label{Hankel}
\phi(\rho,\beta,z)\; =\; \frac{1}{2\pi\ii}\int_H t^{-\beta}\, e^{t + z t^{-\rho}}\, dt,
\end{equation}  
where $H$ is the classical Hankel contour between $e^{-\ii\pi} \infty$ and $e^{\ii\pi}\infty$ encircling the origin counterclockwise. The Wright function was introduced in \cite{W0, W1} for $\rho\ge 0$ in connection with the asymptotic study of partitions, and later studied in \cite{W2} for all $\rho > -1.$ In the recent literature, the function $\phi(\rho, \beta,z)$ with $\rho\in (-1,0)$ is sometimes called the Wright function of the second kind - see \cite{CMP} and the references therein. There is a well-known relationship between the Wright function and the two-parameter Mittag-Leffler function
$$E_{\alpha, \beta} (z) \; =\; \sum_{n\ge 0} \frac{z^n}{\Gamma (\beta + \alpha n)}, \qquad  \beta, z\in \CC,\; \alpha > 0,$$
which is illustrated by the Laplace transforms
$$\int_0^\infty \phi(\rho, \beta, x) \, e^{-xt}\, dx \; =\; E_{\rho, \beta} (1/t),\qquad \beta\in\CC, \; \rho, t > 0$$ 
obtained by direct summation, and
\begin{equation}
\label{MLW}
\int_0^\infty \phi(\rho, \beta, -x) \, e^{-xt}\, dx \; =\; E_{-\rho, \beta-\rho} (-t),\qquad \beta\in\CC, \; \rho\in (-1,0), \;  t > 0
\end{equation}
obtained from Fubini's theorem applied to \eqref{Hankel} and the usual contour representation of $E_{\alpha,\beta}(z).$ We refer to Chapter 18 in \cite{EMOT} and Section 2.3. in \cite{GMP} for more material on this classical relationship. 

The connection between Wright functions and probability theory is relevant for $\rho\in (-1,0].$ In this case indeed, it follows on the one hand from the exponential asymptotic expansion in \cite{W2} that $\phi (\rho, \beta, z)$ is integrable on $\rl^-$ and on the other hand, from \eqref{MLW} and the complete monotonicity of Mittag-Leffler functions on $\rl^-$ originally characterized in \cite{Sch}, that it is positive on $\rl^-$ if $\beta \ge 0$ and $\rho\in (-1,0].$ This shows that the function 
\begin{equation}
\label{Danse}
x\mapsto \varphi_{\a,\b}(x)\; =\; \Gamma (\a +\beta) \, \phi(-\alpha, \beta, -x)
\end{equation}
is the density of a positive random variable $\M_{\a,\beta}$ for all $\a \in [0,1)$ and $\beta \ge 0.$ Expanding the Laplace transform \eqref{MLW} implies that this random variable has Mellin transform 
\begin{equation}
\label{Melle}
\esp [\Mab^s] \; =\; \frac{(1)_s}{(\a+\b)_{\a s}},\qquad s > -1,
\end{equation}
with the usual Pochhammer notation $(a)_s = \Ga(a+s)/\Ga(a)$ for $a > 0$. With the terminology of \cite{Janson}, this shows that $\Mab$ has moments of Gamma type. Observe that \eqref{Melle} can be extended to $\a = 1, \beta > 0$ with $\M_{1,\beta} = \B_{1,\beta}$ where, here and throughout, $\B_{a,b}$ denotes the beta random variable with density
$$\frac{\Ga(a+b)}{\Ga(a)\Ga(b)} \, x^{a-1} (1-x)^{b-1} \Un_{(0,1)} (x).$$
Notice that the right-hand side of \eqref{Melle} is never the Mellin transform of a non-negative function for $\a < 0$ or $\{\a > 0, \beta < 0\}$ since it then vanishes inside the definition strip, or for $\a > 1$ since the support of the function would then be $\{0\}.$ In the following, we will call the set $\{\a\in [0,1],\; \beta \ge 0\}$ the admissible set for $\Mab.$ Except in the degenerate case $\{\a = 1, \beta = 0\}$ with $\M_{1,0} = 1,$ the random variable $\Mab$ has always a smooth density over its support which is given by $\varphi_{\a,\b}(x)$ for $\a\in [0,1)$ and by $\beta (1-x)^{\b -1}\Un_{(0,1)}(x) = \varphi_{1,\beta} (x)$ for $\a = 1.$  It is easy to see from \eqref{Melle} that Supp $\Mab\, =\,\rl^+$ for all $\a\in[0,1)$ and $\b\ge 0.$ Observe also that $\M_{0,\beta} = \L,$ the standard exponential random variable, for all $\b\ge 0.$ When $\a\in (0,1),$ the densities $\varphi_{\a,\b}$ are mostly non-explicit. One purpose of this paper is to show some basic, albeit non-trivial, properties of these densities, which transfer immediately to the Wright functions of the second kind. Our first finding is the following.   

\begin{TheoA} The random variable $\M_{\a,\beta}$ is unimodal.
\end{TheoA} 

This result is easy to prove for $\b\ge\a$ and the density $\varphi_{\a,\b}$ is then non-increasing. For $\b <\a,$ the density is increasing-then-decreasing and the proof, which relies on strong unimodality and a detailed analysis of the distributions with binomial moments introduced in \cite{MP}, is more involved. In the case $\beta = 1 -\alpha,$ the random variable $\M_{\a, 1-\alpha} = \M_\alpha$ has a so-called Mittag-Leffler distribution with moment generating function
\begin{equation}
\label{Géné}
\esp [e^{z \M_\a}] \; =\; \Ea(z),
\end{equation}
where $\Ea(z) = E_{\a, 1} (z)$ stands for the classical Mittag-Leffler function, which we extend to $E_0(z) = 1/(1-z)$ on $\{\Re(z) < 1\}$ for $\a = 0$. If $1-\a \le \a,$ the unimodality of $\Ma$ follows then from the well-known fact - see Section 8 in \cite{Bi2} and the references therein - that it is distributed as the positive part of a real spectrally negative stable random variable with parameter $1/\a\in[1,2],$ which is unimodal by Yamazato's theorem - see Theorem 53.1 in \cite{Sat}. 

The random variable $\M_\a$ plays a central role in various probabilistic contexts. If $\{Z^{(\a)}_t,\; t\ge 0\}$ denotes the $\a-$stable subordinator with Laplace transform $\esp[e^{-\lambda Z^{(\a)}_t}] = e^{-t\lambda^\a}$ for $t,\lambda \ge 0$ and if we set $\Za = Z^{(\a)}_1,$ then it is classical and easy to show by moment identification and self-similarity that
$$\Ma\;\elaw\; \Za^{-\a}\; \elaw\; \inf\{ t > 0, \; Z^{(\a)}_t > 1\}.$$
The Mittag-Leffler random variable also appears in limit theorems for occupation times of Markov processes \cite{Bi1}, P\'olya urn schemes \cite{Jan}, coalescents \cite{Mo}, and elephant random walks \cite{Ber}. This list is non-exhaustive and we refer to the survey paper \cite{Hui} for an account on $\Ma$ and other related random variables. It is natural to ask for basic distributional properties of $\Ma,$ and we focus in this paper on the log-concavity of the density. This property is equivalent to the strong additive unimodality of $\Ma$ in the non-degenerate case $\a\neq 1$ by Ibragimov's theorem \cite{Ibra}, and can hence be viewed as a refinement of unimodality. We prove the following rather unexpected characterization.

\begin{TheoB} The density of $\Ma$ is log-concave if and only if $\a \le \a_*,$ where $\a_*= 0.771667...$ is the unique solution on $(0,1)$ of the equation
$$\frac{1}{\Ga(1-2\a)^2}\; =\; \frac{1}{\Ga(1-\a)\,\Ga(1-3\a)}\cdot$$
\end{TheoB}

Some graphical evidence for the existence and uniqueness of the solution to the above equation defining the critical parameter $\a_*$ is given in Figure \ref{Tuxi} below, where the property is also checked rigorously. For $\a\in (1/2,1),$ the above result has immediate applications to the integro-differential equation
$$u(t,x)\; =\; f(x) \; +\; \frac{1}{\Ga(2\a)}\int_0^t \, (t-s)^{2\a -1}\Delta u(s,x)\, ds,\qquad t > 0,\, x\in\rl,$$ 
which was studied in \cite{Fuj} as an interpolation between the heat equation ($\a = 1/2$) and the wave equation ($\a = 1$), and whose unique solution is expressed in Formula (1) of \cite{Fuji} as
$$u_\a (t,x)\; =\; \frac{1}{2}\lpa \esp \lcr f( x + t^\a \Ma) \, +\, f(x - t^\a \Ma)\rcr\rpa.$$
If the initial data $f(x)$ is unimodal on $\rl,$ as it often happens for natural phenomena governed by heat or waves, then Theorem B and Ibragimov's theorem imply that for every $t > 0,$ the solution $u_\a(t,x)$ will be the superposition of two unimodal functions for $\a\le\a_*,$ and of two functions which may have multiple modes if $\a > \a_*.$ In other words, the spatial behaviour of the solution $u_\a(t,x)$ is similar to the classical solution 
$$u_1 (t,x)\; =\; \frac{1}{2}\lpa \esp \lcr f( x + t) \, +\, f(x - t)\rcr\rpa$$
of the wave equation for $\a\in (1/2,\a_*],$ and possibly different for $\a\in (\a_*,1).$ We refer to \cite{CMP, GMP} and the references therein for more general integro-differential equations and related physical problems, also involving the random variables $\M_{\a,0}, \M_{\a,\a}$ and $\M_{\a,1}.$ 

Our original motivation for Theorem B stems from the so-called generalized entropies on a finite state space. If $p =\{p_i,\; i = 1,\ldots, n\}$ is a probability on a finite set, then its generalized entropy is defined as 
\begin{equation}
\label{Entro}
S(p)\; =\; \sum_{i=1}^n g(p_i)
\end{equation}
where $g$ is some function from $[0,1]$ to $\rl^+$ which satisfies the three Shannon-Khinchin principles. The latter are simply given by: (i) $g(0) = g(1) =0$; (ii) $g$ is continuous; (iii) $g$ is concave, and we refer to \cite{AF, RACF} and the references therein for more on this topic. The classical entropy corresponds to the function $g(x) = - x\log x,$ which involves the natural logarithm. Because of its very simple series representation and its absolute monotonicity on the real line, it is sometimes said - see e.g. Chapters 1 and 2 in \cite{GM} - that the classical Mittag-Leffler function $\Ea(x)$ is for $\a\in [0,1]$ a generalization of the exponential function; we can hence consider its inverse function $\Ea^{-1}(x) = \log_\a(x)$ from $(0,\infty)$ to $\rl$ as a generalized logarithm. For $\a < 1$ however, the polynomial behaviour of $\Ea$ at $-\infty$ makes the continuous extension of $- x\log_\a(x)$ at zero different from the case $\a =1.$ It is easy to see from the first term in the expansion 18.1.(20) in \cite{EMOT} and from the evaluation $\Ea(0) = 1$ that the function
$$g_\a(x) \; =\; -x\log_\a(x) \; +\; \frac{x-1}{\Ga(1-\a)}$$
is well defined on $[0,1]$ with $g_\a(0) = g_\a (1) = 0.$ The continuity of $g_\a$ on $[0,1]$ is plain, and in order to ensure that the functional in \eqref{Entro} with $g = g_\a$ defines a generalized entropy we need to study the concavity of $g_\a,$ which amounts to that of its affine translation $x\mapsto -x\log_\a x.$ 

\begin{TheoC} The function $x\mapsto - x\log_\a(x)$ is concave on $[0,1]$ if and only if $\a\le\a_*$ or $\a = 1.$ 
\end{TheoC}

The reason why the same critical parameter $\a_*$ appears in both Theorem B and C is partly explained by the aforementioned connection between the Wright function $\phi(-\a,1-\a,-x)$ and the classical Mittag-Leffler function $\Ea$. The if part of Theorem C is indeed a natural consequence of the if part of Theorem B, which is the difficult part of the characterization relying on a certain Meijer $G-$function approximation if $\a\le 1/2,$ and on a combination of the so-called Yamazato property and the bell-shape property for real stable densities if $\a > 1/2.$ We refer to Sections 52 and 53 in \cite{Sat} resp. \cite{Kwas, KS} for extended discussions on the Yamazato property resp. the bell-shape property. The fact that the only if parts are the same in both theorems is however more surprising and we comment on this in Remark \ref{Stones} below, where we also discuss in the half-Gaussian case $\a = 1/2$ the equivalence between Theorem C and the classical Sampford inequality on Mill's ratio.  

In view of the previous theorems, it is natural to investigate the log-concavity of the density of $\Mab$ for all admissible $(\a,\b).$ We can show the following partial result. 
  
\begin{TheoD} {\em (a)} The density of $\Mab$ is not log-concave for all $\b\in (0,1)$ and $\a > \a_*(\b),$ where $\a_*(\b)$ is the unique solution on $(0,1)$ of the equation
$$\frac{1}{\Ga(\b-\a)^2}\; =\; \frac{1}{\Ga(\b)\,\Ga(\b-2\a)}\cdot$$

\medskip

\noindent
{\em (b)} The density of $\Mab$ is log-concave in the following situations.

\smallskip

{\em (i)} For all $\a\in [0,4/5] $ and $\b = 0.$ 

\smallskip

{\em (ii)} For all $\a\in [0,1]$ and $\b\ge \a.$ 
 
\end{TheoD}

The necessary condition in Part (a) and the sufficient condition in Part (b) (i), in the non-trivial situation $\a\in (\a_*,4/5],$ are obtained similarly as in Theorem B. However, the sufficient condition given in Part (b) (ii) is more difficult to prove in the relevant case $\a > 4/5.$ The argument relies on the fractional differential rules connecting Wright functions with one another and on a visual bell-shape property strengthening Theorem A, in addition to all tools already appearing in the proofs of Theorem A and B. Unfortunately, it does not seem to us that these tools are enough to tackle the remaining situations $\{0 < \b < \a\}$ and $\{\a > 4/5, \b =0\},$ where we believe that the log-concavity should be characterized by the domain in $(\a,\b)$ appearing in the necessity Part (a). We also conjecture, as in Theorem C, that this property is equivalent to the existence of a generalized entropy associated to the Mittag-Leffler function $E_{\a,\a+\b}$ on the negative half-line. More detail on this open problem is given in the last section of the paper. In this last section we also provide some further results related to the above four theorems, on the multiplicative strong unimodality of $\Mab,$ on the number of positive zeroes of certain Wright functions, and on the reciprocal convexity of Mittag-Leffler functions on the positive half-line. 

\section{Proof of Theorem A} 

The cases $\a = 0,1$ are immediate and we need to consider the situation $\a \in (0,1)$ only. We begin with the easy case $\b\ge \a,$ where we show that $\varphi_{\a,\b}$ is non-increasing and hence unimodal. We use the fractional integration relationship
\begin{equation}
\label{Frac}
\phi(-\a,\beta, -x)\; =\; \frac{1}{\Gamma(\gamma)} \int_x^\infty \phi (-\a, \b - \a\gamma, -t) \, (t-x)^{\gamma -1} \, dt
\end{equation}
which is valid for all $\gamma > 0$ and a straightforward consequence of Fubini's theorem applied to \eqref{Hankel}. Choosing $\gamma = 1$ and differentiating leads to
\begin{equation}
\label{Neg}
\varphi_{\a,\b}'(x)\; =\; -\Gamma(\a+\b)\, \phi(-\a, \b -\a, -x)\; <\; 0
\end{equation}
on $\rl^+$ for $\b\ge \a,$ which completes the argument. We next consider the case $\b < \a.$ Observe that here, one has
$$\varphi_{\a, \b}'(0)\; =\; \frac{(\a -\b)\,\Gamma(\a+\b)}{\Ga(\b)\Ga(\b + 1-\a)}\; >\; 0$$
so that $\varphi_{\a,\b}(x)$ is not non-increasing on $\rl^+$ and we need to show that it increases and then decreases. To prove such a property is usually more difficult than the sole monotonicity and we will need less elementary  tools. In the case $\a\le 1/2$ we use the independent factorization
$$\Mab\;\elaw\; \M_{\a,0} \,\times \, \B_{\a,\b}^\a,$$
which is valid for all admissible parameters as a direct consequence of \eqref{Melle} and a fractional moment identification. The density of $\B_{\a,\b}^\a$ reads
$$\frac{\Ga(\a+\b)}{\Ga(\a + 1)\Ga(\b)}\, (1- x^{1/\a})^{\b -1} \Un_{(0,1)} (x)$$
and is increasing hence unimodal. Besides, setting $f_\a$ for the density of the positive $\a-$stable random variable $\Za,$ one has $\varphi_{\a, 0} (x) = \a x \varphi_{\a, 1-\a} (x) =  x^{-1/\a} f_\a (x^{-1/\a})$ with $f_\a(e^x)$ log-concave on $\rl$ for $\a \le 1/2$ by the main result of \cite{TSPAMS}. Hence, the function $\varphi_{\a, 0}(e^x)$ is log-concave on $\rl$ for $\a \le 1/2$ as well. By the Cuculescu-Theodorescu theorem \cite{CT}, this implies that $\M_{\a, 0}$ is multiplicatively strongly unimodal, in other words the independent product of $\M_{\a, 0}$ with any unimodal random variable remains unimodal. This completes the argument for $\a \le 1/2.$ 

In order to handle the final case $\a\in (1/2,1)$ and $\b\in [0,\a)$ we will need the following lemma, which has an independent interest. Recall that a non-negative function on a given interval is absolutely monotone (AM) if it is smooth and all its derivatives are non-negative.

\begin{Lemm}
\label{Kanter} 
For every $\a\in [1/2,1)$ and $\b\in [0,\a),$ there exists a random variable $\Xab$ such that 
$$\esp[\Xab^s]\; =\;  \frac{(1)_s}{(\a+\b)_{\a s}\,(2-\a-\b)_{(1-\a)s}}$$
for all $s > -1.$ This random variable has a compact support and an {\em AM} density.
\end{Lemm}

This lemma concludes the proof by the same multiplicative strong unimodality argument as above, since \eqref{Melle} gives the independent factorization
\begin{equation}
\label{Fakto}
\Mab\; \elaw\; \G_{2-\a-\b}^{1-\a} \,\times\, \Xab
\end{equation}
where $\Xab$ has an absolutely monotone density and is hence unimodal, whereas $\G_{2-\a-\b}^{1-\a}$ is easily seen to be multiplicatively strongy unimodal by the above criterion. 

\qed

\medskip

\noindent
{\bf Proof of Lemma \ref{Kanter}}. We begin with the case $\a = 1/2.$ First, the Legendre-Gauss multiplication formula for the Gamma function and a fractional moment identification show that the random variable $\X_{1/2,\b}$ does exist and is distributed, with the notation $\B_{a,0} = 1$ for all $a > 0,$ as the independent product 
$$\X_{1/2,\b}\;\elaw\; 2\,\sqrt{\B_{1/2,\b}\,\times\,\B_{1,1/2-\b}},$$
whose support is $[0,2].$ In the case $\b = 0,$ the density is $(x(1-x^2/4)^{-1/2})/2$ and clearly AM. In the case $\b\in(0,1/2),$ the density of the independent product $\B_{1/2,\b}\,\times\,\B_{1,1/2-\b}$ is easily given in terms of the standard hypergeometric function as 
$$\frac{\Ga(\b)\Ga(1/2-\b)}{\pi \,\sqrt{1-x}}\; \pFq{2}{1}{\b,,1-\b}{1/2}{1-x}$$
on $(0,1),$ and the two transformations
\begin{eqnarray*}
\pFq{2}{1}{\b,,1-\b}{1/2}{1-x} & = & \frac{1}{\sqrt{x}}\; \pFq{2}{1}{\b-1/2,,1/2-\b}{1/2}{1-x}\\
& = & \frac{1}{\sqrt{x}}\lpa \sin(\pi\b)\; \pFq{2}{1}{\b-1/2,,1/2-\b}{1/2}{x}\right.\\
& & \qquad\qquad\qquad\qquad\quad\quad \left. \, +\; (1-2\b)\cos(\pi\b)\; \pFq{2}{1}{\b,,1-\b}{3/2}{x}\rpa
\end{eqnarray*}
respectively given by the formulas 2.1.4.(23) and 2.10.(1) in \cite{EMOT}, lead after a change of variable to the following expression for the density of $\X_{1/2,\b}:$
$$\frac{1}{\sqrt{1-z}}\lpa \frac{\Ga(1/2-\b)}{\Ga(1-\b)}\; \pFq{2}{1}{\b-1/2,,1/2-\b}{1/2}{z}\; +\; \frac{(1-2\b)\Ga(\b)}{\Ga(1/2 +\b)}\; \pFq{2}{1}{\b,,1-\b}{3/2}{z}\rpa,$$
with the notation $z = x^2/4.$ The function
$$x\;\mapsto\; \frac{1}{\sqrt{1-z}}\; \pFq{2}{1}{\b,,1-\b}{3/2}{z}\; =\; \pFq{2}{1}{3/2 -\b,,1/2+\b}{3/2}{z}$$
is an entire series in $x^2$ with positive coefficients and hence AM. By Formula 2.1.3.(10) in \cite{EMOT}, so is
$$x\;\mapsto\; \frac{1}{\sqrt{1-z}}\; \pFq{2}{1}{\b-1/2,,1/2-\b}{1/2}{z}\; =\; \frac{\sqrt{\pi}}{\Ga(\b)\Ga(1/2-\b)}\int_0^1 t^{-\b-1/2} (1-t)^{b-1} \frac{(1-tz)^{1/2-\b}}{(1-z)^{1/2}}\, dt$$
because for every $t\in (0,1),$ the function
$$x\;\mapsto\; \frac{(1-tz)^{1/2-\b}}{(1-z)^{1/2}}\; =\; (1-tx^2/4)^{-\b} \sqrt{t + \frac{1-t}{1-x^2/4}}$$
is AM as the product of two AM functions. This shows that the density of $\X_{1/2,\b}$ is AM as the sum of two AM functions.

We next consider the more involved case $\a\in (1/2,1).$ The existence of the random variable $\Xab$ for all $\b\in [0,\a)$ follows from the main result of \cite{MP} which implies that there exists an absolutely continuous and  compactly supported random variable $\Y_{\a,r}$ with Mellin transform
$$\esp[ \Y_{\a,r}^s] \; = \; \frac{(r+1)_s}{(1)_{\a s} \, (r+1)_{(1-\a) s}},\qquad s > -r$$
for all $\a\in (0,1)$ and $r\in (-1,1/\a-1].$ Setting $r = (1-\b)\a^{-1}-1$ shows that $\Xab \elaw \Y_{\a,r}^{(-r)}$ has the required Mellin transform, and a compact support $[0,\a^{-\a}(1-\a)^{\a-1}].$

It remains to show that the density of the random variable $\Xab$ is AM. We first suppose $\a = l/k$ rational with $k,l$ positive integers, $k\ge 3$ and $l < k < 2l.$ Applying Theorem 3.1 in \cite{MP} and (3.8) therein shows the independent product representation
$$\Y_{\a,r}\; \elaw\; \a^{-\a}(1-\a)^{\a -1}\, \lpa \B_{\b_1, \tilde{\a}_1 - \b_1}\,\times\,\cdots\,\times\,\B_{\b_k, \tilde{\a}_k - \b_k}\rpa^{1/k}$$
with the notation of \cite{MP} for the sequence $\{(\b_j, \tilde{\a}_j), j = 1,\ldots, k\},$ that is
$$\b_j\, =\, \frac{r+j}{k} \qquad\mbox{and}\qquad \tilde{\a}_j\, =\, \lacc\begin{array}{ll}\frac{i}{l} & \mbox{if $j = j'_i, 1\le i\le l,$} \\
\frac{r+j-i}{k-l} & \mbox{if $j'_i < j < j'_{i+1}$}
\end{array}\right.\quad\mbox{with}\quad j'_i\, =\, \left\lfloor \frac{ik}{l} -r\right\rfloor.$$
It is not difficult to see from the proof of Lemma 3.2 in \cite{MP} that one has $\b_j\le\tilde{\a}_j\le\b_{j+1}$ for all $j =1,\ldots, k$ so that this sequence is interlacing, that is $\beta_1\,\le \tilde{\a}_1 \le \beta_2\le\ldots\le \beta_k\le \tilde{\a}_k.$ Since we have $\b_j -\b_i = (j-i)/k\not\in\ZZ$ for all $i\neq j,$ the density of the product $\B_{\b_1, \tilde{\a}_1 - \b_1}\times\ldots\times\B_{\b_k, \tilde{\a}_k - \b_k}$ is a Meijer $G$-function having the following convergent series representation 
\begin{equation}
\label{Meije}
\lpa\prod_{i=1}^k\frac{\Ga(\tilde{\a}_i)}{\Ga(\b_i)}\rpa \sum_{j=1}^k \frac{x^{\b_j-1}}{\Ga(\tilde{\a}_j - \b_j)}\prod_{i\neq j} \frac{\Ga(\b_i - \b_j)}{\Ga(\tilde{\a}_i - \b_j)}\sum_{n= 0}^\infty \lpa \prod_{m=1}^k \frac{(1 +\b_j - \tilde{\a}_m)_n}{(1+\b_j-\b_m)_n}\rpa x^n
\end{equation}  
on $(0,1)$ - see Formula (2.1) in \cite{Dunk}. For all $i\neq j\in\{1,\ldots, k\},$ the interlacing property shows that either $\b_i -\b_j > 0$ and $\tilde{\a}_i - \b_j > 0$ or $\b_i - \b_j \in (-1,0)$ and $\tilde{\a}_i - \b_j \in (-1,0],$ so that 
$$\prod_{i\neq j} \frac{\Ga(\b_i - \b_j)}{\Ga(\tilde{\a}_i - \b_j)}\; \ge\; 0$$
for all $j = 1,\ldots, k.$ Moreover, for all $j,m\in \{1,\ldots, k\}$ we have $1 +\b_j -\b_m\ge 1+\b_1 - \b_k = 1/k > 0$ on the one hand, and 
$$1 + \b_j - \tilde{\a}_m \, \ge \, 1 +\b_1 - \tilde{\a}_k \, =\, \frac{1+r}{k}\, - \,\lpa\frac{r}{k-l}\rpa_+$$ 
on the other hand, where the equality comes from (3.6) in \cite{MP} and the right-hand side is non-negative since either $r\in (-1, 0],$ or $r > 0$ with
$$\frac{1+r}{k}\, - \,\lpa\frac{r}{k-l}\rpa_+ =\; \frac{k-l(1+r)}{k(k-l)}\; =\; \frac{\b}{k-l}\; \ge \; 0.$$ 
This implies that all coefficients in \eqref{Meije} are non-negative. Putting everything together, we have shown that the density of $\Xab \elaw \Y_{\a,r}^{(-r)}$ is given by the convergent series representation
$$\sum_{j=1}^k \sum_{n=0}^\infty a_{j,n}\, x^{k(n+\b_j) -1 -r}\; =\; \sum_{j=1}^k \sum_{n=0}^\infty a_{j,n}\, x^{kn +j -1}$$
for some $a_{j,n}\ge 0,$ where in the equality we have used, with the notation of \cite{MP},
$$k(n+\b_j) - 1 -r\, = \, k(n+ \b_j - \b_1) + k\b_1 - 1-r\, = \, k(n+ \b_j -\b_1) \, =\, kn + j-1$$
for all $n\ge 0$ and $j \in\{1,\ldots, k\}.$ This shows the AM character of the density of $\X_{\a,\b}$ and concludes the proof for $\a$ rational. The general case follows from a standard approximation argument using Mellin inversion. We omit details.   

\qed
 
\begin{Rem} 
\label{Pense}
{\em (a) In the case $\b = 1-\a,$ the random variable $\X_{\a,1-\a}$ has a very simple Mellin transform given by 
$$\esp[\X_{\a,1-\a}^s]\; =\;  \frac{(1)_s}{(1)_{\a s}\,(1)_{(1-\a)s}}$$
for all $s > -1,$ which is invariant by the switching $\a\leftrightarrow 1-\a.$ Several properties of this random variable, which is an explicit deterministic transform of the uniform random variable, are derived in Section 3 of \cite{TSEJP}. See in particular Proposition 3.1 resp. Remark 2 (b) therein for the increasing character of the density, resp. for a question which was left open on the convexity of the density. The above Lemma \ref{Kanter} shows the more general AM property.
     
\medskip

(b) Setting $f_{\a,\b}$ for the density of $\Xab$, the multiplicative convolution formula and a change of variable give 
$$\frac{\Ga(\a+\b)}{\Ga(\b)} \; =\; \varphi_{\a,\b}(0) \; =\; \frac{1}{\Ga(2-\a)}\int_0^\infty f_{\a,\b} (0) \, e^{-y}\,dy,$$
which leads to the following formula, which can also be retrieved from the results of \cite{MP}:
$$f_{\a,\b} (0)\; =\; \frac{\Ga(\a +\b)\, \Ga(2-\a)}{\Ga(\b)}\cdot$$
The behaviour of $f_{\a,\b}$ at the right-hand boundary $\a^{-\a}(1-\a)^{\a-1}$ of its support can be obtained from the first term of Formula (1.4) in \cite{Dunk} in the case when $\a$ is rational. Skipping details, we get
$$f_{\a,\b}(x)\; \sim\; \frac{(1-\a)^{2\a +\b -1} \,(1-\a-\b)}{\sqrt{2}\, \a^\b\, \sin (\pi(\a +\b))}\, (1-\a^\a(1-\a)^{1-\a} x)^{-1/2}$$
for $\b\neq 1-\a$ and 
$$f_{\a,1-\a}(x)\; \sim\; \frac{(1-\a)^{\a}}{\sqrt{2}\,\pi\, \a^{1-\a}}\, (1-\a^\a(1-\a)^{1-\a} x)^{-1/2},$$
as $x\to \a^{-\a}(1-\a)^{\a-1}.$ Since the constants do not depend on $k,l$ these behaviours easily extend to the case when $\a$ is not rational. A curious feature is the "universality" of $-1/2$ as a power exponent, which comes from the fact, with the notations of \cite{Dunk} and \cite{MP}, that
$$\delta \; =\; \sum_{i=1}^k \tilde{\a}_i - b_i\; =\; \sum_{i=1}^k \a_i -\b_i\; =\; 1/2$$
for all values of $\a = l/k\in[1/2,1)$ rational and $\b\in [0,\a).$}
\end{Rem}
 
\section{Proof of Theorem B} 

We begin with the easy only if part. By \eqref{Danse}, we have 
$$(\varphi_{\a,1-\a}'(0))^2\; -\; \varphi_{\a,1-\a}(0)\,\varphi_{\a,1-\a}''(0)\; =\; \frac{1}{\Ga(1-2\a)^2}\, -\, \frac{1}{\Ga(1-\a)\,\Ga(1-3\a)}\, = \, \rho(\a)$$
and $\Ma$ is not log-concave if $\rho(\a)$ is negative. The function $\rho(\a)$ is clearly non-negative for $\a\in (0,2/3],$ as the sum of two non-negative functions for $\a\in [1/3,2/3]$ and by the log-convexity of the Gamma function for $\a \in (0,1/3).$ We hence need to show the existence of a unique $\a_* \in (2/3,1)$ such that $\rho(\a) < 0$ for $\a \in (\a_*,1)$ and $\rho(\a) > 0$ for $\a\in (2/3,\a_*).$ Factorizing
$$\rho(\a)\; =\; \frac{1}{\Ga(1-2\a)^2}\,\lpa 1\, -\,\frac{(3\a -1) (3\a -2)\,\Ga(2(1-\a))^2}{(2\a-1)^2\,\Ga(1-\a)\,\Ga(3(1-\a))}\rpa,$$
we see that the function
$$\a\,\mapsto\,\frac{(3\a -1) (3\a -2)\,\Ga(2(1-\a))^2}{(2\a-1)^2\,\Ga(1-\a)\,\Ga(3(1-\a))}\; =\; \lpa 9\, -\, \frac{1}{4(2\a-1)^2}\rpa\times\lpa\frac{\Ga(2(1-\a))^2}{4\,\Ga(1-\a)\,\Ga(3(1-\a))}\rpa$$
increases on $(2/3,1)$ from $0$ to $3$ by the log-convexity of the Gamma function. This finishes the proof of the only if part.

\medskip

\begin{figure}
\centering

\includegraphics[scale =0.9]{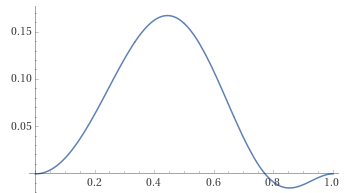}
\caption{Plot of $\a\mapsto \rho(\a)$ for $\a\in (0,1).$}
\label{Tuxi}

\end{figure}

We now proceed to the proof of the if part. The argument is different according as $\a\le 1/2$ or $\a > 1/2.$ In the former case it hinges upon the a.s. convergent infinite product representation
\begin{equation}
\label{Mainar}
\Ma\;\elaw\; \frac{1}{\Ga(1+\a)}\,\prod_{n=0}^\infty \lpa \frac{n+1}{n+\a}\rpa \B_{1+\frac{n}{\a}, \frac{1}{\a}-1}
\end{equation}
which is a consequence of the main result in \cite{JSW} in the case $m=1$ and equation (3) therein. This implies that we are reduced to show the log-concavity on $(0,1)$ of the density of the finite independent product
$$\Z_{n,\a}\; =\; \B_{1, \frac{1}{\a} -1}\,\times\,\B_{1+\frac{1}{\a}, \frac{1}{\a}-1}\, \times\, \cdots\, \,\times\, \B_{1+\frac{n}{\a}, \frac{1}{\a}-1}$$
for all $n.$ Fix $n\ge 0$ and let $\psi_{k,n}$ be the density of $\B_{1+\frac{n-k}{\a}, \frac{1}{\a}-1}\, \times\, \cdots\, \,\times\, \B_{1+\frac{n}{\a}, \frac{1}{\a}-1}$ for all $k =0,\ldots, n.$ We prove by induction on $k$ that the function $\tilde{\psi}_{k,n}(x) = x^{\frac{k-n}{\a}} \psi_{k,n}(x)$ is log-concave for all $k=0,\ldots, n,$ which will complete the argument since the density of $\Z_{n,\a}$ is $\tilde{\psi}_{n,n}.$ The case $k = 0$ is straightforward because
$$\tilde{\psi}_{0,n}(x) \; =\; \frac{\Ga(\frac{n+1}{\a})}{\Ga(\frac{1}{\a} -1) \Ga (1+ \frac{n}{\a})}\; (1-x)^{\frac{1}{\a} -2} \, \Un_{(0,1)} (x),$$
which is log-concave since $\a \le 1/2.$ Supposing next $\tilde{\psi}_{k,n}$ log-concave for some $0\le k < n$, the multiplication convolution
$$\psi_{k+1,n} (x)\; =\; \int_x^1 \psi_{k,n} (y)\, \psi_{0,n-k-1} (xy^{-1})\, \frac{dy}{y}$$
implies 
$$\tilde{\psi}_{k+1,n} (x)\; =\; \frac{\Ga(\frac{n-k}{\a})}{\Ga(\frac{1}{\a} -1) \Ga (1+ \frac{n-k-1}{\a})}\;\int_\rl y \,\tilde{\psi}_{k,n} (y)\, (y-x)^{\frac{1}{\a} -2}\Un_{(-1,0)}(x-y) \, dy,$$
which is log-concave by the Prékopa-Leindler theorem since both functions  $(-y)^{\frac{1}{\a} -2}\Un_{(-1,0)}(y)$ and $y \tilde{\psi}_{k,n} (y)$ are log-concave on $\rl,$ the former by $\a\le 1/2$ and the latter by the induction hypothesis.

\medskip

For $\a > 1/2,$ we will take advantage of the aforementioned fact that $\Ma$ is the positive part of a real spectrally positive $(1/\a)-$stable random variable. We will use the fundamental and recently established property that the density $\varphi_\a$ of this random variable is bell-shaped, that is the $n-$th derivative $\varphi_\a^{(n)}$ vanishes exactly $n$ times on $\rl$ for every $n\ge 1$ - see Corollary 1.3 in \cite{Kwas}. It is well-known that 
$$\int_0^\infty \varphi_\a (x)\, dx\; =\; \a$$
so that $\varphi_\a (x) = \a \varphi_{\a, 1-\a} (x)$ for all $x\ge 0.$ We hence need to show that 
$$\psi_\a (x)\; =\; \varphi_\a'(x)^2\, -\, \varphi_\a(x)\varphi_\a''(x)\; \ge\; 0$$
for all $\a\in(1/2,\a_*]$ and $x\ge 0.$ We will distinguish three cases. First, we suppose $\a\in (1/2, 2/3]$ with
$$\varphi_\a'(0)\; =\; \frac{-\a}{\Ga(1-2\a)}\; >\; 0\qquad\mbox{and}\qquad \varphi_\a''(0)\; =\; \frac{\a}{\Ga(1-3\a)}\; \le\; 0.$$
The bell-shape property implies that there exists $0 < a_1 < b_2$ such that $(a_1-x)\varphi_\a'(x) \ge 0$ and $(x-b_2)\varphi_\a''(x) \ge 0$ for all $x\ge 0.$ This implies on the one hand that $\psi_\a(x) \ge 0$ for all $x\in [0,b_2].$ On the other hand, the Yamazato property applied to the real spectrally positive stable density $\varphi_\a(-x)$ with unique mode $-a_1$ shows that this function is log-concave on $(-\infty,a_1].$ Indeed, with the notation of \cite{SY} we have $k(x) = c_\a x^{-1-1/\a}\Un_{\{x>0\}}$ so that $\lambda_- = 0, \lambda_+ = \infty$ and $\varphi_\a(-x)$ is of the type $I_7,$ and we can apply Theorem 1.3.(xii) therein. This shows that $\psi_\a(x) \ge 0$ for all $x\in [a_1,\infty)$ as well and completes the proof since $a_1 < b_2.$ 

We next suppose $\a\in (2/3, 3/4]$ with $\varphi_\a'(0) > 0,$ 
$$\varphi_\a''(0)\; =\; \frac{\a}{\Ga(1-3\a)}\; >\; 0\qquad\mbox{and}\qquad \varphi_\a'''(0)\; =\; \frac{-\a}{\Ga(1-4\a)}\; \le\; 0.$$
The bell-shape property implies that there exists $0 < b_1 < a_1 < b_2$ and $0 < b_1 < c_2 < b_2 < c_3$ such that $(b_1-x)(b_2 -x)\varphi_\a''(x) \ge 0$ and $(x-c_2)(c_3 -x)\varphi_\a'''(x)\ge 0$ for all $x\ge 0.$ In particular $\varphi_\a$ is concave on $[b_1,b_2]$ and $\psi_\a(x) \ge 0$ for all $x\in[b_1,b_2].$ Recalling that $\varphi_\a$ is log-concave on $[a_1,\infty)$ with $a_1 < b_2$ by the Yamazato property, we are reduced to show that $\psi_\a(x)\ge 0$ for all $x\in [0,b_1].$ This follows from
$$\psi_\a (0)\; =\; \frac{\a^2}{\Ga(1-2\a)^2}\, -\, \frac{\a^2}{\Ga(1-\a)\Ga(1-3\a)}\; > \; 0$$
where the inequality is a consequence of $\a\le 3/4 < \a_*,$ and from
$$\psi_\a' (x)\; =\; \varphi_\a'(x)\varphi_\a''(x)\, -\, \varphi_\a(x)\varphi_\a'''(x)\; \ge\; 0\qquad\mbox{for all $x\in[0,b_1]$}$$
because $\varphi_\a,\varphi_\a'$ and $\varphi_\a''$ are non-negative on $[0,b_1]$ whereas $\varphi_\a'''$ is non-positive on $[0,c_2].$ 

We finally suppose $\a\in (3/4,\a_*]$ with $\varphi_\a'(0) > 0, \varphi_\a''(0) >0,$ 
$$\varphi_\a'''(0)\; =\; \frac{-\a}{\Ga(1-4\a)}\; >\; 0\qquad\mbox{and}\qquad \varphi_\a''''(0)\; =\; \frac{\a}{\Ga(1-5\a)}\; < \; 0,$$
where the last inequality comes from $3/5 < \a\le \a_* < 4/5.$ The bell-shape property implies that there exists $0 < b_1 < a_1 < b_2$ and $0 < c_1 < b_1 < c_2 < b_2 < c_3$ such that $(b_1-x)(b_2 -x)\varphi_\a''(x) \ge 0$ and $(c_1-x)(c_2-x)(c_3 -x)\varphi_\a'''(x)\ge 0$ for all $x\ge 0.$ As above, we need to show that $\psi_\a(x)\ge 0$ for all $x\in [0,b_1]$ only. The condition $\a\le\a_*$ yields $\psi_\a(0)\ge 0$ and we have $\psi_\a'(x)\ge 0$ for all $x\in[c_1,b_1]$ as above, so that we are reduced to show that $\psi_\a'(x)\ge 0$ for all $x\in[0,c_1]$ as well. On the one hand, the bell-shape property shows that there exist $d_1 < 0 < c_1 < d_2 < c_2 < d_3 < c_3 < d_4$ such that $(d_1-x)(d_2-x)(d_3 -x)(d_4-x)\varphi_\a''''(x)\ge 0$ for all $x\ge 0.$ In particular, we have $\varphi_\a''''\le 0$ on $[0,c_1]$ and
$$\psi_\a''(x)\; =\; (\varphi_\a''(x))^2\, -\, \varphi_\a(x)\varphi_\a''''(x)\; \ge\; 0\qquad\mbox{for all $x\in[0,c_1].$}$$  
Finally, we compute
\begin{eqnarray*}
\psi_\a' (0) & = & \frac{\a^2}{\Ga(1-\a)\Ga(1-4\a)}\, -\, \frac{\a^2}{\Ga(1-2\a)\Ga(1-3\a)}\\
& = & \, -\, \frac{2\a^2}{\Ga(1-2\a)\Ga(1-3\a)}\lpa \frac{1}{2}\, -\, \frac{(4\a-1)(4\a-3)\,\Ga(2(1-\a))\,\Ga(3(1-\a))}{(3\a-1)(3\a-2)\,\Ga(1-\a)\,\Ga(4(1-\a))}\rpa\; >\; 0,
\end{eqnarray*}
where the last inequality comes from 
$$\frac{(4\a-1)(4\a-3)\,\Ga(2(1-\a))\,\Ga(3(1-\a))}{(3\a-1)(3\a-2)\,\Ga(1-\a)\,\Ga(4(1-\a))}\, =\, \lpa 16 \, -\,\frac{5}{(3\a-1)(3\a-2)}\rpa\, \times\,\frac{\Ga(2(1-\a))\,\Ga(3(1-\a))}{9\,\Ga(1-\a)\,\Ga(4(1-\a))}$$
which increases as a function of $\a\in[3/4,1)$ from 0 to 1, taking the value
$$\frac{11\, \Ga(2/5)\Ga(3/5)}{14\,\Ga(1/5)\,\Ga(4/5)}\; =\; \frac{11}{7(1+\sqrt{5})}\; <\; 1/2$$
at $\a =4/5 >\a_*.$ This completes the proof.

\qed

\begin{Rem}
{\em A direct proof of the log-concavity of $\Ma$ can be given in the explicit or semi-explicit cases $\a = 1/3, 1/2, 2/3.$ For $\a = 1/2,$ the argument is immediate since $\varphi_{1/2,1/2} (x) = e^{-x^2/4}/\sqrt{\pi}.$ For $\a = 1/3,$ the second order ODE
$$3\,\varphi_{1/3,2/3}'' (x) \; =\; x\, \varphi_{1/3,2/3} (x)$$
implies that the function $\psi_{1/3}(x) = \varphi_{1/3,2/3}'(x)^2 - \varphi_{1/3,2/3} (x) \varphi_{1/3,2/3}''(x)$  is such that 
$$\psi_{1/3}'(x)\, = \,-\varphi_{1/3,2/3}(x)^2/3\, <\, 0\quad\mbox{with}\quad \left\{\begin{array}{l} \psi_{1/3}(0)\, =\, 1/\Ga(1/3)^2\, >\, 0\\
\lim_{x\to\infty} \psi_{1/3}(x)\, =\, 0,\end{array}\right.$$ 
so that $\psi_{1/3}(x) > 0$ for all $x \ge 0$ as required. Alternatively, we can use the representation $\varphi_{1/3,2/3} (x) = 3^{2/3} {\rm Ai} (3^{-1/3} x)$ - see e.g. Section 3.2. in \cite{CMP}, and the fact that $x\mapsto {\rm Ai} (x)$ is log-concave on $(-a_1,\infty)$ where $-a_1 < 0$ is its first negative zero - see Proposition 2 in \cite{Salm}. In the case $\a = 2/3,$ we have the Whittaker representation
$$\varphi_{2/3,1/3}(x)\; =\; 2^{-1/3} \sqrt{\frac{3}{\pi}} \,e^{-X/2} X^{-1/3}\, \cW_{1/2,1/6} (X)\; =\; 2^{-1/3} \sqrt{\frac{3}{\pi}}\, e^{-X} X^{1/3}\, \Psi(1/6,4/3,X)$$
where $X = 4x^3/27$ and
$$\Psi (a, c, z) \; = \; \frac{1}{\Ga(a)}\,\int_0^\infty e^{-zt} t^{a-1} (1+t)^{c-a-1} dt,\qquad a,z > 0, c\in\rl,$$
is a confluent hypergeometric function - see e.g. Section 3.3 in \cite{CMP} and formulas 6.5.(2) and 6.9.(2) in \cite{EMOT}. Setting $\psi(x) =  \Psi(1/6,4/3,x)$ for concision, we see after some elementary simplifications that the log-concavity of $\varphi_{2/3,1/3}(x)$ amounts to 
$$(6x+1)\psi(x)^2\, +\, 9((\psi'(x)^2 - \psi(x) \psi''(x))\; \ge\; 6x\, \psi(x)\psi'(x), \qquad x\ge 0.$$ 
By the confluent hypergeometric equation $6 x \psi''(x) = (6x - 8) \psi'(x) + \psi(x),$ this is equivalent to
$$2 \, (\psi(x) + 3x\, \psi'(x))^2\, +\, 9 x\, \psi(x)^2\; \ge\; 18 x^2\,\psi(x)\psi'(x), \qquad x\ge 0,$$
which holds true because the left-hand side is positive and the right-hand side is negative. }
\end{Rem}

\section{Proof of Theorem C}

We discard the classical case $\a = 1$ with $(x \log x)'' = 1/x > 0$ on $(0,1).$ The case $\a = 0$ is also immediate with $\log_0(x) = 1 -1/x$ and $(x \log_0 x)'' = 0.$ For $\a\in (0,1),$ we have 
$$(x \log_\a x)''\; = 2\log_\a'(x)\, + \,x \log_\a''(x)\; =\; \frac{2}{\Ea'(\log_\a(x))}\, -\,\frac{x\,\Ea''(\log_\a(x))}{\Ea'(\log_\a(x))^3}\cdot$$
Changing the variable $x =\Ea(z)$ with $z\in\rl^-$ and using $\Ea'(z) > 0$ for all $z\in\rl^-,$ we deduce the equivalence
\begin{equation}
\label{ConvConv}
-x \log_\a x\;\;\mbox{is concave on $(0,1)$}\quad\Leftrightarrow\quad \Ea\Ea''\, \le\, 2 (\Ea')^2\;\; \mbox{on $\rl^-$}\quad\Leftrightarrow\quad \frac{1}{\Ea} \;\;\mbox{is convex on $\rl^-.$}
\end{equation}
The asymptotic expansion 18.1.(20) in \cite{EMOT} of the Mittag-Leffler function for $z\to -\infty$ yields
$$\Ea(z)\; =\; -\sum_{n=1}^3 \frac{z^{-n}}{\Ga(1-\a n)}\, +\, O(\vert z\vert^{-4}), \qquad \Ea'(z)\; =\; \frac{1}{\a} \, E_{\a,\a} (z)\; =\; \sum_{n=1}^3 \frac{n\, z^{-n-1}}{\Ga(1-\a n)}\, +\, O(\vert z\vert^{-5})$$
and
$$\Ea''(z)\; =\; \frac{1}{\a^2} \lpa E_{\a,2\a-1} (z)\, +\, (1-z) E_{\a, 2\a}(z)\rpa\; =\; -\sum_{n=1}^3 \frac{n(n+1)\, z^{-n-2}}{\Ga(1-\a n)}\, +\, O(\vert z\vert^{-6}).$$
After simplification, this leads to
$$\lim_{z\to -\infty} \, z^6\lpa 2 \,\Ea'(z)^2\, -\, \Ea(z)\Ea''(z)\rpa\; =\; 2\lpa \frac{1}{\Ga(1-2\a)^2} \, -\, \frac{1}{\Ga(1-\a)\Ga(1-3\a)}\rpa$$
and by the proof of Theorem B, the right-hand side is negative for $\a\in(\a_*,1).$ By the first equivalence in \eqref{ConvConv}, this concludes the proof of the only if part.

We now proceed to the proof of the if part, which will be a consequence of Theorem B. By the second equivalence in \eqref{ConvConv}, we need to show that the function
$$x\;\mapsto\; \frac{\Ea'(-x)}{\Ea(-x)^2}$$
is non-decreasing on $\rl^+.$ On the one hand, we have
$$\Ea'(-x)\; =\; \int_0^\infty e^{-xt} \,t\, \varphi_{\a,1-\a} (t)\, dt.$$
On the other hand, the additive convolution formula and a change of variable entail
$$\Ea(-x)^2 \; =\; \int_0^\infty e^{-xt} \, t\lpa\int_0^1 \varphi_{\a,1-\a} (tu)\,\varphi_{\a,1-\a} (t(1-u))\, du\rpa dt.$$
Hence, applying Lemma 3 in our former paper \cite{FS} - see also the references therein, we are reduced to show that the function
$$t\;\mapsto\; \int_0^1 \lpa \frac{\varphi_{\a,1-\a} (tu)\,\varphi_{\a,1-\a} (t(1-u))}{\varphi_{\a,1-\a} (t)}\rpa du$$ 
is non-decreasing on $\rl^+.$ But since $\a\le\a_*,$ it is a straightforward consequence of the log-concavity of $\varphi_{\a,1-\a}$ that the function
$$t\;\mapsto\; \frac{\varphi_{\a,1-\a} (tu)\,\varphi_{\a,1-\a} (t(1-u))}{\varphi_{\a,1-\a} (t)}$$
is non-decreasing for every $u\in [0,1].$ This completes the proof.

\qed

\begin{Rem} 

\label{Stones}

{\em (a) A combination of Theorem B, Theorem C and \eqref{ConvConv} implies the following equivalence for every $\a\in (0,1):$
\begin{equation}
\label{Surprise}
\a\,\le\,\a_*\;\Longleftrightarrow\; \Ma\;\,\mbox{has a log-concave density}\;\Longleftrightarrow\; \frac{1}{\Ea(-t)}\; \,\mbox{is convex on $\rl^+.$}
\end{equation}
The second equivalence is surprising in view of Corollary 1.2 in \cite{KL}, which shows that the log-concavity of $\Ma$ is actually equivalent to the convexity on $\rl^+$ of all functions $\Ea^{(n)}\!(-t)^{-1/n}, \, n\ge 1.$ In particular, we retrieve that for all $\a\le\a_*$ the function
$$\frac{1}{\Ga(\a)\, E_{\a,\a}(-t)}\; =\; \frac{1}{\Ga(\a + 1)\, \Ea'(-t)}\; = \; \sqrt{\frac{1}{\Ga(\a + 1)\, \Ea'(-t)}}\times \sqrt{\frac{1}{\Ga(\a + 1)\, \Ea'(-t)}}$$
is convex as the product of two increasing convex functions. This is actually true for all $\a\le 4/5$ by the log-concavity of $\varphi_{\a, 0}$ established in Theorem D and the above reasoning. Another interesting feature of \eqref{Surprise} is that the condition $\a\le\a_*$ is read off from the sole behaviour of $\varphi_{\a,1-\a}(x)$ at zero resp. of $\Ea(-t)$ at $\infty.$ One can check that this is also governed by the behaviour at infinity of each function $\Ea^{(n)}\!(-t)^{1/n}, \, n\ge 1.$

\medskip

(b) It is worth recalling from \eqref{Géné} that by Hölder's inequality, the function $\Ea(x)$ is log-convex on $\rl$ for all $\a \in (0,1)$ as the moment generating function of the random variable $\Ma$. Hence, the inequality $\Ea'(x)^2 \le \Ea(x) \Ea''(x)$ holds for all $x\in \rl$ and $\a\in (0,1).$ The equivalence \eqref{ConvConv} and Theorem C show that
$$\Ea'(x)^2 \;\le \;\Ea(x) \Ea''(x)\; \le\; 2\Ea'(x)^2$$ 
holds for all $x\in\rl^-$ and $\a\le\a_*.$ In Paragraph \ref{Conv} below, we will see that it also holds for all $x\in\rl^+$ and $\a\in (0,1).$

\medskip

(c) As a consequence of Theorem C, the function  
\begin{equation}
\label{CM}
x\;\mapsto\; \frac{\Ea(-x^\a)^2}{\Ea'(-x^\a)}
\end{equation}
is non-increasing on $\rl^+$ if and only if $\a\le \a_*.$ In our previous paper \cite{FS}, we proved that it is also completely monotone for $\a = 1/2,$ as a consequence of the hypergeometric identity in law
$$\frac{\G_{1/2}}{\G_{1/2}}\; +\; \frac{\G_{1/2}}{\G_{1/2}}\; \elaw\; \frac{\G_1}{\G_{1/2}}\; \times\;\lpa 1\, +\, \sqrt{\B_{1/2,1/2}}\rpa.$$ 
Indeed, we have
$$E_{1/2} (-x) \; =\; e^{x^2} \,{\rm Erfc} (x)\; =\; \sqrt{\frac{\pi}{2}}\; r(\sqrt{2} x)$$
where $r$ stands for the so-called Mill's ratio, and we can apply Proposition 10 in \cite{FS} which is a consequence of the former identity. In this respect, Theorem C for $\a = 1/2$ turns out to be equivalent to the classical Sampford inequality for Mill's ratio - see Proposition 11 in \cite{FS}. Observe that Proposition 10 in \cite{FS} also shows that the function in \eqref{CM} is logarithmically completely monotone for $\a = 1/2.$ We believe that this property remains true for all $\a\le 1/2.$}
\end{Rem}

\section{Proof of Theorem D} 

\subsection{Proof of Part (a)} This is similar to the proof of the only if part of Theorem B. By \eqref{Danse}, we have 
$$\varphi_{\a,\b}'(0)^2\, -\, \varphi_{\a,\b}(0)\varphi_{\a,\b}''(0)\; =\; \Ga(\a+\b)^2 \lpa \frac{1}{\Ga(\b-\a)^2}\, -\, \frac{1}{\Ga(\b)\,\Ga(\b-2\a)}\rpa$$
and $\Mab$ is not log-concave if the function on the right-hand side is negative. This function is positive for $\a\le(\b + 1)/2$ and converges to $\b^2(\b-1) < 0$ as $\a \to 1.$  Hence, we need to ckeck that the function
$$\a\,\mapsto\, \frac{\Ga(\b)\,\Ga(\b-2\a)}{\Ga(\b-\a)^2}\; =\; \lpa \frac{(\a-\b)^2}{(2\a -\b) (2\a -\b-1)}\rpa\times\lpa \frac{\Ga(\b)\,\Ga(\b+2(1-\a))}{\Ga(\b +(1-\a))^2} \rpa$$
decreases on $((\b+1)/2,1].$ The second factor on the right-hand side is positive and decreasing by the log-convexity of the Gamma function, and the same is true for the first factor because its logarithmic derivative equals
$$\frac{\b(\a-\b) -\a(1-\b)}{(\a-\b)(2\a -\b)(2\a -\b-1)}\; < \; 0.$$

\qed

\subsection{Proof of Part (b)} 

\subsubsection{Case {\em (i)}} First, we observe that the case $\a\le \a_*$ is an immediate consequence of Theorem B because $\phi(-\a, 0, -t) = \a t \,\phi(-\a, 1-\a, -t)$ is log-concave if $\phi(-\a,1-\a,-t)$ is log-concave. We are hence reduced to show that $\M_{\a,0}$ is log-concave for every $\a\in (\a_*, 4/5].$  Setting $\tva(x) = x\varphi_\a(x)$ and 
$$\tpa (x)\; = \;\tva'(x)^2 \, - \, \tva(x)\tva''(x)\; = \; x^2\,\psi_\a(x)\, + \,\varphi_\a(x)^2$$ 
with the notation of the proof of Theorem B, we need to show that $\tpa (x) \ge 0$ for every  $x\ge 0.$ 
The Yamazato property for $\Ma,$ which is valid for every $\a\in [1/2,1),$ implies $\tpa(x)\ge \varphi_\a(x)^2 > 0$ for all $x\in [b_1,\infty)$ and we need to show that $\tpa(x) \ge 0$ for all $x\in[0,b_1].$ We have 
$$\tpa'(x)\; =\; 2\varphi_\a(x)\varphi_\a'(x)\, +\, 2x\, \psi_\a(x)\, +\, x^2\psi_\a'(x)\; \ge\; 2\varphi_\a(x)\varphi_\a'(x)\, +\, 2x\, \psi_\a(x)$$
for all $x\in [0,b_1]$ since we have seen during the proof of Theorem B that $\psi_\a'(x) \ge 0$ on this interval. Setting finally $\tppa(x) = 2\varphi_\a(x)\varphi_\a'(x)\, +\, 2x\, \psi_\a(x),$ we have $\tppa(0) = 2\varphi_\a(0)\varphi_\a'(0) > 0$ since $\a > 1/2,$ and
$$\tppa'(x)\; = \; 4\varphi_\a'(x)^2\; +\; 2x\, \psi_\a'(x)\; >\; 0$$
for all $x\in [0,b_1].$  This implies that $\tpa'(x) > 0$ for all $x\in [0,b_1]$ and concludes the proof because $\tpa (0) = \varphi_\a(0)^2 > 0.$

\qed

\begin{Rem} 
\label{DD}
{\em (a) Unfortunately, the above argument does not convey directly to $\a\in(4/5,1)$ since we then need to consider the further derivatives of $\tpa(x)$ and the variational study becomes increasingly complicated as $\a$ approaches 1. Let us mention a discretization argument which is similar to the proof of the case $\a\le 1/2$ in Theorem B, and which works as well for $\a\le 4/5.$ It relies on the factorization
$$\M_{\a,0}\;\elaw\; \frac{\Ga(\a)}{\Ga(2\a)}\,\prod_{n=0}^\infty \lpa \frac{n+1+\a}{n+2\a}\rpa \B_{2+\frac{n}{\a}, \frac{1}{\a}-1}$$
which is a consequence of \eqref{Mainar} and $\M_{\a,0}\elaw\Ma^{(1)}$ with the standard notation for size-bias. Setting
$$\W_{n,\a}\;\elaw\; \B_{2, \frac{1}{\a} -1}\,\times\,\B_{2+\frac{1}{\a}, \frac{1}{\a}-1}\, \times\, \cdots\, \,\times\, \B_{2+\frac{n}{\a}, \frac{1}{\a}-1}$$
for all $n\ge 0,$ we are then reduced to show that $\W_{4p +3,\a}$ is log-concave for all $p\ge 0$ and $\a\le 4/5.$ With the notation of the proof of Theorem B, this is done in computing the density of $\W_{3,\a}$ up to some constant as $x\,(1-x)^{\frac{4}{\a} -5} \, (1+x)\,(1 + c_\a x + x^2)\,\Un_{(0,1)} (x)$ with 
$$c_\a\; = \; \frac{10(1-\a)}{3-\a}$$
on the one hand, and in showing that the bivariate function $(x,y)\mapsto x (x+y)(x^2 + c_\a xy + y^2)$
is log-concave on $\{0\le x\le y\}$ if and only if $c_\a \ge 10/11\Leftrightarrow\a\le 4/5$ on the other hand. In general, the density of $\Z_{n,\a}$ exhibits a curious family of polynomials having positive and symmetric coefficients, which will be the matter of further research.\\

(b) The fact that $\varphi_{\a, 0}$ is log-concave and $\varphi_{\a, 1-\a}$ is not log-concave for $\a\in (\a_*,4/5]$ show that fractional integration of order in $(0,1)$ does not preserve log-concavity in general, in view of
$$\varphi_{\a, 1-\a}(x) \; =\; \frac{1}{\Ga(\a)\,\Ga(\gamma)} \int_x^\infty \varphi_{\a, 0}(t) \, (t-x)^{\gamma -1}\, dt$$
with $\gamma = \a^{-1} - 1 \in (0,1).$ We are not aware of other counterexamples in the literature. } 

\end{Rem}

\subsubsection{Case {\em (ii)}} We first observe that the case $\a \le 4/5$ is a direct consequence of Case (i) and the Prékopa-Leindler theorem in view of 
\begin{equation}
\label{a0}
\varphi_{\a,\a}(x)\; =\; \frac{\Ga(2\a)}{\Ga(\a)}\int_x^\infty \varphi_{\a,0} (t)\, dt.
\end{equation}
We begin with the case $\b = \a > 4/5.$ Setting $a > 0$ for the unique mode of $\varphi_{\a,0}$ on $(0,\infty),$ which exists by Theorem A, we have
$$\varphi_{\a,\a}'(x)^2\, -\, \varphi_{\a,\a}(x) \varphi_{\a,\a}''(x)\; =\; \varphi_{\a,0}(x)^2\, +\, \frac{\Ga(2\a)}{\Ga(\a)}\,\varphi_{\a,\a}(x) \varphi_{\a,0}'(x),$$
which is clearly positive on $[0,a].$ Moreover, the Yamazato property for $\Ma,$ which is valid since $\a > 1/2,$ shows that 
$$\varphi_{\a,0} (x)\; = \;\Ga(\a +1)\, x\, \varphi_{\a,1-\a} (x)$$ 
is log-concave on $[b,\infty)$ where $b< a$ is the unique mode of $\varphi_{\a,1-\a},$ and the Prékopa-Leindler theorem transfers this property to $\varphi_{\a,\a}$ by \eqref{a0}. This shows that $\varphi_{\a,\a}$ is log-concave on $[0,a]\cup[b,\infty) = \rl^+$ and concludes the argument for the boundary case $\b =\a.$ Clearly, this implies also the property for $\b \ge 2\a$ by the relationship
$$\varphi_{\a,\b}(x)\; =\; \frac{\Ga(\a +\b)}{\Ga(2 \a)\,\Ga(\gamma)}\,\int_x^\infty \varphi_{\a,\a} (t)\, (t-x)^{\gamma - 1} \,dt$$
with $\gamma = \b\a^{-1} - 1\ge 1.$  

We next handle the case $\b \in (\a, 1]$ with $\a > 4/5.$ The formula 
\begin{equation}
\label{second}
\varphi_{\a,\b}'' (x)\; =\; -\frac{\Ga(\a+\b)}{\Ga(\b)}\, \varphi_{\a,\b-\a}' (x),
\end{equation}
implies by Theorem A that $\varphi_{\a,\b}$ is concave and hence log-concave on $[0,c]$ where $c$ is the unique mode of $\varphi_{\a,\b-\a},$ which is positive since $\b-\a < \a.$ Moreover, for $\b\in (\a,1)$ we have 
$$\varphi_{\a,1-\a}'(x)\; =\; \frac{1}{\Ga(\b)\Ga(\gamma)} \int_x^\infty \varphi_{\a, \b - \a}'(t) \, (t-x)^{\gamma -1} \, dt$$
with $\gamma = (1-\beta)/\a,$ which readily implies by positivity that $c > b$ since $b$ is the unique mode of $\varphi_{\a,1-\a}(x).$ The Yamazato property for $\Ma$ and the formula 
$$\varphi_{\a,\b}(x)\; =\; \frac{1}{\Ga(\gamma)} \int_x^\infty \varphi_{\a, 0}(t) \, (t-x)^{\gamma -1} \, dt$$
with $\gamma =\b\a^{-1} > 1$ imply by the Prékopa-Leindler theorem that $\varphi_{\a,\b}(x)$ is log-concave on $[b, \infty)$ as well, and hence on $[0,c]\cup[b,\infty)= \rl^+,$ which concludes the proof for $\b\in(\a,1),$ and the boundary case $\b = 1$ follows by continuity.

We finally consider the case $\b\in (1,2\a)$ with $\a > 4/5,$ which is unfortunately technical. The difficulty stems from the fact that we cannot use neither a combination of \eqref{Frac} and the lower boundary case $\b = 1$ because of Remark \eqref{DD} (b), nor the simple modality argument of the case $\b\in(\a,1]$ since then $c\to 0$ as $\b\to 2\a.$ Instead, the proof will rely on the following visual refinement of Theorem A, which has an independent interest. We recall that an inflection point of a smooth real function is a point where the second derivative has a strict change of sign. 

\begin{Lemm}
\label{Kanter2}
For every $\a\in [1/2,1]$ and $\b\in [0,\a],$ the density of $\M_{\a,\b}$ has at most two inflection points.
\end{Lemm}
  
Postponing the proof of this Lemma to the end of the section, we first terminate the proof. Recall that since $\b > \a,$ we have $\varphi_{\a,\b}' < 0$ by \eqref{Neg}. Moreover, by \eqref{second} and Theorem A, there exists $a > 0$ such that $(x-a)\varphi_{\a,\b}''(x)\ge 0.$ Setting
$$\psi_{\a,\b}(x)\; =\; \varphi_{\a,\b}'(x)^2\; -\; \varphi_{\a,\b}(x) \varphi_{\a,\b}''(x)$$
as in the proof of Theorem B, we have $\psi_{\a,\b}(x)\ge 0$ for $x\in [0,a].$ Moreover, Lemma \eqref{Kanter2} and \eqref{second} imply that the function $\varphi_{\a,\b}'''(x)$ vanishes at most twice on $\rl^+$ and hence only once on $(a,\infty).$ This clearly shows that there exists $b > a$ such that $\varphi_{\a,\b}'''(x)\ge 0$ for $x\in[a,b)$ and $\varphi_{\a,\b}'''(x)\le 0$ for $x\in[b,\infty),$ and we obtain that
$$\psi_{\a,\b}'(x)\; =\; \varphi_{\a,\b}'(x)\varphi_{\a,\b}''(x)\; -\; \varphi_{\a,\b}(x) \varphi_{\a,\b}'''(x)$$
is non-positive on $[a,b].$ All in all, we are hence reduced to show that $\psi_{\a,\b}$ is non-negative on $[b,\infty).$ Since
$$\varphi_{\a,\b}(x)\; =\; \frac{\Ga(\a+\b)}{\Ga(\gamma)}\int_x^\infty \varphi_{\a,1-\a}(t)\, (t-x)^{\gamma -1} dt$$
with $\gamma = 1 + (\b-1)\a^{-1} > 1,$ the Prékopa-Leindler theorem and the Yamazato property for $\Ma$ show that $\psi_{\a,\b} (x) \ge 0$ for all $x\in [b_1,\infty)$ where $b_1$ is the first inflection point of $\varphi_{\a,1-\a}.$ It remains to show that $b_1\le b$ and for this we use
$$\varphi_{\a,1-\a}'(x)\; =\; \frac{1}{\Ga(\gamma)\,\Ga(\a+\b)} \int_x^\infty \varphi_{\a,\b}'''(t)\, (t-x)^{\gamma -1}\, dt$$
with $\gamma = 1 + (1-\b)\a^{-1} > 0,$ which yields $\varphi_{\a,1-\a}'(x) \le 0$ for $x\ge b$ and hence $b_1 < a_1 \le b$ as required, where $a_1$ is the unique zero of $\varphi_{\a,1-\a}'$. 
 
\qed
 
\medskip

\noindent
{\bf Proof of Lemma \ref{Kanter2}}. For the simplicity of notation, we rewrite the factorization \eqref{Fakto} as 
$$\Mab\; \elaw\; \a^{-\a}(1-\a)^{\a-1}\,\G_{2-\a}^{1-\a} \,\times\, \Yab$$
with $\Xab = \a^{-\a}(1-\a)^{\a-1}\Yab$ and the random variable $\Yab$ having support $[0,1]$ and an absolutely monotone density on $(0,1)$ by Lemma \ref{Kanter}. Taking e.g. a partition of unity, for every $\varepsilon \in (0,1)$ we consider the approximation $\Yabe$ having a smooth density on $\rl^+$ vanishing on $[1,\infty),$ equalling that of $\Yab$ on $[0,1-\varepsilon],$ and joining $1-\varepsilon$ to $1$ with one mode and two inflection points. 

The AM character of the density of $\Yab$ implies the crucial property that the density $f_{\a,\b, \varepsilon}$ of $\Yabe$ has one mode and two inflection points on the whole $[0,1].$ Setting now
$$\Mabe\; \elaw\; \G_{2-\a}^{1-\a} \,\times\, \Yabe,$$
it is enough by approximation to show that the density $\varphi_{\a,\b, \varepsilon}$ of $\Mabe$ has at most two inflection points for every $\varepsilon > 0.$ By multiplicative convolution, we have
$$\varphi_{\a,\b, \varepsilon}(x) \; =\; \int_0^1 f_{\a,\b, \varepsilon} (y)\, g_\a(xy^{-1})\, \frac{dy}{y}$$
where $g_\a(y)$ stands for the density of $\G_{2-\a}^{1-\a}.$ The exponential behaviour of $g_\a$ at infinity allows one to differentiate inside the integral: we get 
$$\varphi_{\a,\b, \varepsilon}'(x) \; =\; \int_0^1 f_{\a,\b, \varepsilon} (y)\, g_\a'(xy^{-1})\, \frac{dy}{y^2}\; =\; \int_0^1 f_{\a,\b, \varepsilon}' (y)\, h_\a (xy^{-1})\, \frac{dy}{y}$$
with the notation $h_\a(t) = t^{-1} g_\a(t),$ where the second equality follows from an integration by parts obtained from the vanishing character of $f_{\a,\b,\varepsilon}(y)$ at $y=1$ resp. of $g_\a(xy^{-1})$ at $y = 0.$ Differentiating again and using an analogous integration by parts, we obtain
$$\varphi_{\a,\b, \varepsilon}''(x) \; =\; \int_0^1 f_{\a,\b, \varepsilon}' (y)\, h_\a'(xy^{-1})\, \frac{dy}{y^2}\; =\; \frac{1}{x} \int_0^1 f_{\a,\b, \varepsilon}'' (y)\, h_\a (xy^{-1})\, dy.$$
Now since
$$h_\a(t) = \frac{t^{\frac{\a}{1-\a}}\, e^{-t^{\frac{1}{1-\a}}}}
{(1-\a)\,\Ga(2-\a)},$$
it follows from the basic example (2.1) p. 15 in \cite{Karl} that the kernel $h_\a(xy^{-1})$ is strictly totally positive on $(0,\infty)\times (0,\infty).$ And since $\sharp\{ x > 0,\; f_{\a,\b, \varepsilon}''(x) = 0\}\, =\, 2$ by construction, the variation-diminishing property given in Theorem 3.1.(ii) of \cite{Karl} shows that $\varphi_{\a,\b, \varepsilon}''(x)$ vanishes at most twice on $(0,\infty)$ for every $\varepsilon > 0,$ which completes the proof.

\qed

\section{Further remarks}
\label{Furth}

\subsection{Multiplicative strong unimodality} The proof of Theorem A depends on the notion of multiplicative strong unimodality (MSU), introduced in \cite{CT}. For a positive random variable $X$ with density $f$, it is easy to check that the MSU property is implied by the non-increasing character of $f$, in other words by the unimodality at zero. When $f$ is not non-increasing, the MSU property is equivalent by the main result of \cite{CT} to the log-concavity of $t\mapsto f(e^t)$ and this property, which we may call MSU*, is invariant under power transformations of $X$. Observe that a positive random variable with non-increasing density needs not be MSU*. We refer to \cite{CT} and the introduction of \cite{TSPAMS} for more material on multiplicative strong unimodality. The following characterization, which is a simple consequence of Theorem D, improves on the main result of \cite{TSPAMS}.

\begin{Propo}
\label{MSU}
For every $(\a,\b)$ admissible, one has
$$\Mab\,\mbox{is {\em MSU}}\;\Longleftrightarrow\;\Mab\,\mbox{is {\em MSU*}}\;\Longleftrightarrow\;\b\ge \a\,\;\mbox{or}\,\; \{\b = 0, \a\le 1/2\}.$$
\end{Propo}

\proof Suppose first $0 < \b < \a.$ Then we have $\varphi_{\a,\b}'(0) > 0$ and $\varphi_{\a,\b}$ is not non-increasing on $\rl^+$ so that MSU and MSU* are equivalent properties, tantamount to the non-increasing character of
$$x\;\mapsto\; \frac{x\,\varphi_{\a,\b}'(x)}{\varphi_{\a,\b} (x)}\cdot$$
Differentiating, the properties are equivalent to $x\varphi_{\a,\b}(x)\varphi_{\a,\b}''(x) + \varphi_{\a,\b}(x)\varphi_{\a,\b}'(x)\, \le\, x \varphi_{\a,\b}'(x)^2,$
which is not true in the neighbourhood of zero since $\varphi_{\a,\b}(0)\varphi_{\a,\b}'(0) > 0.$ Supposing next $\b\ge\a,$ we have seen that $\varphi_{\a,\b}'\le 0$ on $\rl^+$ so that $\M_{\a,\b}$ is MSU. Moreover, it is also MSU* since
$$x\varphi_{\a,\b}(x)\varphi_{\a,\b}''(x)\, + \,\varphi_{\a,\b}(x)\varphi_{\a,\b}'(x)\; \le\; x\varphi_{\a,\b}(x)\varphi_{\a,\b}''(x)\;\le\; x \varphi_{\a,\b}'(x)^2,$$
where in the second inequality we have used Theorem D (b) (ii). Finally, if $\b =0,$ then 
$\varphi_{\a,0}(0) = 0$ and MSU and MSU* are again equivalent properties. But since 
$$\varphi_{\a,0} (e^t)\; =\; \Ga(1+\a)\, e^t\, \varphi_{\a,1-\a} (e^t),$$ 
the MSU* property of $\M_{\a,0}$ amounts to that of $\Ma,$ and we can apply the main result of \cite{TSPAMS}. 

\endproof

\begin{Rem} 
\label{Duran}
{\em Another direct consequence of Theorem D (b) (ii) and the differential rule \eqref{Frac} is the following Tur\'an inequality for the Wright function, which holds true for all $\a\in(0,1)$ and $\b,x\ge 0:$ 
$$\phi(-\a,\b+\a, -x)\,\phi(-\a, \b-\a,-x)\;\le\; \phi(-\a, \b,-x)^2.$$
Observe that the latter inequality is equivalent to the log-concavity of the sequence 
$$\{\phi(-\a, \b +(n-1)\a,-x),\; n\ge 0\}$$ 
for all $\a\in(0,1)$ and $\b, x\ge 0.$ There is a vast recent literature on Tur\'an inequalities for special functions and we refer to \cite{Cso} for a survey. See also \cite{KK} and the references therein for results and conjectures in the framework of Meijer $G-$functions. Notice that in many papers on Tur\'an inequalities, the coefficients in the involved series have constant sign, which is not the case in our framework by the negativities of $-\a$ and $-x$.}
\end{Rem}

\subsection{On the number of positive zeroes of the Wright function} In this paragraph, we briefly comment on how some of our findings can be applied to investigate the exact number of positive zeroes of Wright functions, in a non-probabilistic setting. Our results are only partial as a simple consequence of Theorem A and Lemma \ref{Kanter2}, and this topic certainly deserves a deeper study. For any $\rho > 0,$ the Wright function $\phi(\rho,\b,z)$ has an infinite number of zeroes given by the Hadamard factorization, and it was shown in Theorem 1 of \cite{ABS} that the zeroes are all negative for $\b > 0$ and that there are $[\b] +1$ non-negative zeroes for $\b < 0$ and $\rho\in(0,1].$ The situation is however more complicated for $\rho\in(-1,0),$ since then the Wright function has order greater than one - see \cite{W2}. Setting $\cN_{\rho,\b,+} = \sharp\{ x > 0, \; \phi(\rho,\b,-x)\, =\, 0\},$ we see by the discussion made in the introduction that 
$$\cN_{\rho,\b,+}\; =\; 0\quad\Longleftrightarrow\quad \b\,\ge\, 0$$ 
for every $\rho\in (-1,0).$ In the non-probabilistic setting, we can show the following.

\begin{Propo}
\label{Zero}
For every $\rho\in (-1,-1/2]$, one has
$$\cN_{\rho,\b,+}\; =\; 1\quad\mbox{for $\b\in[-1,0)$}\qquad\mbox{and}\qquad\cN_{\rho,\b,+}\; =\; 2\quad\mbox{for $\b\in[2\rho,-1).$}$$
\end{Propo}
 
\proof If $\rho\in (-1,0)$ and $\b\in[\rho,0),$ we know from \eqref{Frac} that 
$$\phi(\rho,\b, -x)\; =\; -\phi'(\rho,\b-\rho,-x)$$
and the function on the right-hand side is negative-then-positive on $(0,\infty)$ by Theorem A, with an isolated zero on $(0,\infty)$ by the real analyticity of Wright functions. If $\rho\in (-1,-1/2]$ and $\b\in(-1,\rho),$ we know from \eqref{Frac} that 
\begin{equation}
\label{Second}
\phi(\rho,\b, -x)\; =\; \phi''(\rho,\b-2\rho,-x)
\end{equation}
and the function on the right-hand side is also negative-then-positive on $(0,\infty),$ by Lemma \ref{Kanter2} and the fact that $\phi(\rho,\b,0) = 1/\Ga(\b) < 0,$ with again an isolated zero on $(0,\infty).$ If $\rho\in (-1,-1/2]$ and $\b = -1,$ then $\phi(\rho,-1,0) = 0$ and $\phi(\rho,-1, -x) < 0$ in a neighbourhood of zero, so that again $\cN_{\rho,-1,+} =1.$ Finally, if $\rho\in (-1,-1/2)$ and $\b\in[2\rho,-1),$ then \eqref{Second} and Lemma \ref{Kanter2} with $\phi(\rho,\b,0) = 1/\Ga(\b) > 0$ show that $\phi(\rho,\b,-x)$ is positive-then-negative-then-positive on $(0,\infty),$ with two isolated zeroes.

\endproof

\begin{Rem} 
\label{Positano}
{\em In the non-probabilistic setting and for $\rho\in(-1,0),$ there exist Wright functions which never vanish on the positive half-line. For $\rho\in(-1,-1/2]$ and $\b > 1+\rho,$ the function 
$$\phi (\rho,\b,x)\; =\; \frac{1}{\Ga(\gamma)} \int_{-\infty}^x \phi(\rho,1+\rho,-t) \,(x-t)^{\gamma-1}\, dt$$
with $\gamma = 1 + (1-\b)\rho^{-1} >0,$ never vanishes on $(0,\infty)$ since $\phi(\rho, 1+\rho,-t)$ is the density of a spectrally negative stable random variable with index $-1/\rho$, which is positive on $\rl.$ However  the function is not integrable at infinity by the asymptotic behaviour
$$\phi (\rho,\b,x)\; \sim\;- \frac{x^{(1-\b)/\rho}}{\rho\,\Ga(1+ (1-\b)/\rho)}\qquad\mbox{as $x\to\infty,$}$$
with $(1-\b)/\rho > -1,$ which is a consequence of Theorem 4 in \cite{W2}. Observe also from the same asymptotics that $\phi(\rho,\b,x)$ is integrable and takes negative values on $(0,\infty)$ for $\rho\in(-1,-1/2]$ and $\b < 1+\rho.$ We refer to \cite{Louche} for an asymptotic study of the zeroes of the Wright function, relying on the original asymptotic behaviours given in \cite{W1,W2}.}

\end{Rem}

\subsection{On the reciprocal convexity of the Mittag-Leffler function}
\label{Conv}
 As discussed in the above Remark \ref{Stones} (i), for every $\a\in (0,1)$ the reciprocal convexity of the Mittag-Leffler function $E_\a(x)$ on the negative half-line is characterized by $\a\le\a_*,$ as a by-product of Theorem C. From the asymptotic expansion 18.1.(20) in \cite{EMOT} for $\a\in (1,2)$ and the works of Wiman \cite{Wim} for $\a\ge 2,$ it is known that the function $E_\a(x)$ vanishes at least once, and hence cannot be reciprocally convex, on $\rl^-$ for all $\a > 1.$ We hence have the characterization
$$\frac{1}{E_\a}\;\,\mbox{is convex on $\rl^-$}\qquad \Longleftrightarrow\qquad \a\,\le\a_*\quad\mbox{or}\quad \alpha = 1.$$
It is interesting to mention that the reciprocal convexity property holds true on the positive half-line for all two-parameter Mittag-Leffler functions $E_{\a,\b}$ with $\a,\b\ge 0.$ Overall, the argument is much simpler than for the if part of Theorem C by the positivity of the coefficients. 
 
\begin{Propo} 
\label{Convex}
For every $\a,\b \ge 0,$ the function 
$$x\;\mapsto \; \frac{1}{E_{\a,\b}(x)}$$
is convex on $\rl^+.$
\end{Propo}

\proof The case $\a = 0$ with $1/E_{\a,\b}(x) =\Ga(\b) (1-x)_+$ is immediate and we will suppose henceforth $\a > 0.$ We need to prove that the function
$$x\;\mapsto\; \frac{E_{\a,\b} (x)^2}{E_{\a,\b}' (x)}$$
is non-decreasing on $\rl^+.$ Since the coefficients in the entire series 
$$E_{\a,\b} (x)^2\; =\; A(x)\: =\; \sum_{n\ge 0} a_n\, x^n\qquad\mbox{and}\qquad E_{\a,\b}'(x) \; =\; B(x)\; =\; \sum_{n\ge 0} b_n\, x^n$$ are non-negative, it is enough to show by a classical lemma of \cite{BK} that the sequence $\{u_n,\; n\ge 0\}$ is non-decreasing, where
$$u_n\; = \;\frac{a_n}{b_n}\; = \; \frac{\Ga(\b + \a (n +1))}{n+1}\, \sum_{k=0}^n \,\frac{1}{\Ga(\b +\a k)\, \Ga(\b +\a(n-k))}\cdot$$ 
Indeed, this condition entails that the coefficients of the series $A'(x) B(x) - A(x)B'(x)$ are then all non-negative. For every $n\ge 0,$ we compute
\begin{eqnarray*}
u_n - u_{n-1} \!& = &\! \frac{1}{n(n+1)}\lpa  \sum_{k=0}^n \frac{n\, \Ga(\b + \a (n +1))}{\Ga(\b +\a k)\, \Ga(\b +\a(n-k))}\, -\, \sum_{k=0}^{n-1} \frac{(n+1)\,\Ga(\b + \a n)}{\Ga(\b +\a k)\, \Ga(\b +\a(n-1-k))} \rpa\\
& = &\! \frac{2}{n(n+1)}\,\sum_{k=1}^n \,\frac{k}{\Ga(\b +\a(n-k))} \lpa \frac{\Ga(\b + \a (n+1))}{\Ga(\b +\a k)}\, -\, \frac{\Ga(\b + \a n)}{\Ga(\b +\a (k-1))} \rpa\;\ge\; 0,
\end{eqnarray*} 
where the second equality follows from the identities
$$\sum_{k=0}^n\, n\, v_k v_{n-k}\; =\; 2 \;\sum_{k=1}^n\, k\, v_k v_{n-k}\qquad\mbox{and}\qquad \sum_{k=0}^{n-1} \, (n+1)\, v_k v_{n-1-k}\; =\; 2 \;\sum_{k=1}^n \, k\, v_{k-1} v_{n-k}$$
which are valid for any sequence $\{v_n,\; n\ge 0\},$ and the inequality comes from the log-convexity of the Gamma function. This completes the proof.

\endproof

\begin{Rem}
\label{Ger}
{\em For $\a\ge 1$ and $\b\le\a,$ it is easy to see that the sequence 
$$\frac{(n+1)\, \Ga(\b +\a n)}{\Ga(\b+\a +\a n)}$$
is non-increasing and so the function $E_{\a,\b}(x)$ is log-concave on $\rl^+$ by the same argument as above, which implies the reciprocal convexity of $E_{\a,\b}(x)$ on $\rl^+.$ On the other hand, observe that the function $E_{\a,\b}(x) = \Ga(\b) \esp[e^{x\M_{\a, \b-\a}}]$ is log-convex on $\rl^+$ for $\a\le 1$ and $\b\ge\a.$ We refer to \cite{GS} for further log-concavity, log-convexity, subadditivity and superadditivity properties of the classical Mittag-Leffler function on the positive half-line.}  
\end{Rem}

\subsection{A conjecture} The necessary condition given in Theorem D (a) shows that the function $\phi(-\a,\b, -x)$ is not log-concave on $\rl^+$ for $\a > \a_*(\b).$ It is easy to show that this can only happen in the subset $\{0 < \b < 2\a -1 \le 1\}$ of the admissible set, and that it is equivalent to the condition
$$\frac{\Ga(\b)\Ga(\b-2\a)}{\Ga(\b-\a)^2}\; <\; 1$$
inside this subset. Moreover, from the asymptotics 18.1.(20) in \cite{EMOT} we have
$$\lim_{z\to -\infty} \, z^6\lpa 2 \,E_{\a,\a+\b}'(z)^2\, -\, E_{\a,\a+\b}(z)E_{\a,\a+\b}''(z)\rpa\; =\; 2\lpa \frac{1}{\Ga(\b-\a)^2} \, -\, \frac{1}{\Ga(\b)\Ga(\b-2\a)}\rpa$$
and this shows that the condition $\a\le\a_*(\b)$ is also necessary for the reciprocal convexity of $E_{\a,\b}$ on $\rl^-.$ Setting $\a_*(0) = \a_*(\b) = 1$ for all $\b\ge 1,$ in view of \eqref{Surprise} it is very natural to formulate the following.

\begin{Conj} For all admissible $(\a,\b),$ one has
$$\a\,\le\,\a_*(\b)\;\Longleftrightarrow\; \Mab\;\,\mbox{has a log-concave density}\;\Longleftrightarrow\; \frac{1}{E_{\a,\a+\b}(-t)}\; \,\mbox{is convex on $\rl^+.$}$$
\end{Conj}
 
From the beginning of the proof of Theorem C, we also see that the third property amounts to the concavity on $[0,1]$ of the function
$$g_{\a,\b}(x) \; =\; -x\log_{\a,\b}(x)\; +\; \frac{x-1}{\Ga(\b)\Ga(\a+\b)}$$
where $\log_{\a,\b}(x) = E_{\a,\a+\b}^{-1}(x)/\Ga(\a+\b)$ is a generalized logarithm, in other words to the property that the functional in \eqref{Entro} defined with $g_{\a,\b}$ instead of $g_\a$ is a generalized entropy. Theorems B, C and D show the veracity of this conjecture in the cases $\b\ge\a, \b =1-\a$ and $\{\a\le 4/5, \b = 0\}.$ To handle the remaining cases $0 <\b <\a$ and $\{\a > 4/5, \b = 0\}$ however, it seems to the authors that other arguments are needed. In particular, in the case $\b > 1-\a$ we believe, in spite of Remark \ref{Positano}, that $\Mab$ is the positive part of some real self-decomposable spectrally negative random variable, which should obey both the Yamazato property and the bell-shape. We leave this question open for future research. 
  
\section*{Acknowledgement}
Rui A. C. Ferreira was supported by the ``Funda\c{c}\~{a}o para a Ci\^encia e a Tecnologia (FCT)" through the program ``Stimulus of Scientific Employment, Individual Support-2017 Call" with reference CEECIND/00640/2017.

\end{document}